\documentclass[11pt]{amsart}
\usepackage{amsmath}
\usepackage{amsxtra}
\usepackage{amscd}
\usepackage{amsthm}
\usepackage{amsfonts}
\usepackage{amssymb}
\usepackage{eucal}

\newtheorem{lem}{Lemma}

\theoremstyle{definition}

\newtheorem{thm}{Theorem}

\theoremstyle{remark}

\begin{document}

\renewcommand{\baselinestretch}{0.833}

\newcommand{\thmref}[1]{Theorem~\ref{#1}}
\newcommand{\secref}[1]{Sect.~\ref{#1}}
\newcommand{\lemref}[1]{Lemma~\ref{#1}}
\newcommand{\propref}[1]{Proposition~\ref{#1}}
\newcommand{\corref}[1]{Corollary~\ref{#1}}
\newcommand{\remref}[1]{Remark~\ref{#1}}
\newcommand{\conjref}[1]{Conjecture~\ref{#1}}
\newcommand{\nc}{\newcommand}
\nc{\on}{\operatorname}
\nc{\ch}{\mbox{ch}}
\nc{\Z}{{\mathbb Z}}
\nc{\C}{{\mathbb C}}
\nc{\pone}{{\mathbb C}{\mathbb P}^1}
\nc{\pa}{\partial}
\nc{\F}{{\mathcal F}}
\nc{\arr}{\rightarrow}
\nc{\larr}{\longrightarrow}
\nc{\al}{\alpha}
\nc{\ri}{\rangle}
\nc{\lef}{\langle}
\nc{\W}{{\mathcal W}}
\nc{\la}{\lambda}
\nc{\ep}{\epsilon}
\nc{\su}{\widehat{{\mathfrak s}{\mathfrak l}}_2}
\nc{\sw}{{\mathfrak s}{\mathfrak l}}
\nc{\g}{{\mathfrak g}}
\nc{\h}{{\mathfrak h}}
\nc{\n}{{\mathfrak n}}
\nc{\N}{\widehat{\n}}
\nc{\G}{\widehat{\g}}
\nc{\De}{\Delta_+}
\nc{\gt}{\widetilde{\g}}
\nc{\Ga}{\Gamma}
\nc{\one}{{\mathbf 1}}
\nc{\z}{{\mathfrak Z}}
\nc{\zz}{{\mathcal Z}}
\nc{\Hh}{{\mathcal H}_\beta}
\nc{\qp}{q^{\frac{k}{2}}}
\nc{\qm}{q^{-\frac{k}{2}}}
\nc{\La}{\Lambda}
\nc{\wt}{\widetilde}
\nc{\wh}{\widehat}
\nc{\qn}{\frac{[m]_q^2}{[2m]_q}}
\nc{\cri}{_{\on{cr}}}
\nc{\kk}{h^\vee}
\nc{\sun}{\widehat{\sw}_N}
\nc{\hh}{{\mathbf H}_{q,t}}
\nc{\HH}{{\mathcal H}_{q,t}}
\nc{\hhh}{{\mathcal H}_{q,1}}
\nc{\ca}{\wt{{\mathcal A}}_{h,k}(\sw_2)}
\nc{\si}{\sigma}
\nc{\gl}{\widehat{{\mathfrak g}{\mathfrak l}}_2}
\nc{\el}{\ell}
\nc{\s}{T}
\nc{\bib}{\bibitem}
\nc{\om}{\omega}
\nc{\WW}{\W_\beta}
\nc{\scr}{{\mathbf S}}
\nc{\ab}{{\mathbf a}}
\nc{\rr}{r}
\nc{\ol}{\overline}
\nc{\con}{qt^{-1} + q^{-1}t}
\nc{\den}{q^{\el-1} t^{-\el+1}+ q^{-\el+1} t^{\el-1}}
\nc{\ds}{\displaystyle}
\nc{\B}{B}
\nc{\A}{A^{(2)}_{2\el}}
\nc{\GG}{{\mathcal G}}
\nc{\UU}{{\mathcal U}}
\nc{\MM}{{\mathcal M}}
\nc{\CC}{{\mathcal C}}
\nc{\GL}{^L\G}
\nc{\gL}{^L\g}
\nc{\dzz}{\frac{dz}{z}}
\nc{\Res}{\on{Res}}
\nc{\rep}{{\mathcal R}ep \;}
\nc{\uqg}{U_q \G}
\nc{\uqgg}{U_q \g}
\nc{\mc}{\mathcal}
\nc{\Cal}{\mathcal}
\nc{\calp}{{\Cal P}}
\nc{\bp}{{\mathbf P}}
\nc{\bq}{{\mathbf Q}}
\nc{\bb}{{\mathfrak b}}
\nc{\uqb}{U_q \bb_-}
\nc{\uqn}{U_q \wt{{\mathfrak n}}}
\nc{\uqh}{U_q \wt{{\mathfrak h}}}
\nc{\uqhh}{U_q \wh{{\mathfrak h}}}
\nc{\uqnn}{U_q \wh{{\mathfrak n}}}
\nc{\ot}{\otimes}
\nc{\R}{{\mc R}}
\nc{\uqbb}{U_q \wt{\g}}
\nc{\yy}{{\mc Y}}
\nc{\uqsl}{U_q \widehat{\sw}_2}
\nc{\ga}{\gamma}
\nc{\Ab}{{\mathbf A}}
\nc{\Yb}{{\mathbf Y}}
\nc{\yb}{{\mathbf y}}
\nc{\uh}{U \wt{{\mathfrak h}}}
\nc{\uhh}{U \wh{{\mathfrak h}}}
\newcommand{\scs}{\scriptstyle}
\nc{\us}{\underset}
\nc{\opl}{\oplus}
\nc{\yyy}{\wh{\yy}}
\nc{\ovl}{\overline}
\nc{\beq}{\begin{equation}}
\nc{\Fq}{{\mathcal F}}
\nc{\Mq}{{\mathcal M}}
\nc{\Rep}{\on{Rep}}

\nc{\ghat}{\G}
\nc{\mb}{\mathbf}

\nc{\ka}{\kappa}


\title{Affine Kac-Moody algebras, integrable systems and their
deformations}\thanks{The text of the Hermann Weyl Prize lecture given
at the XXIV Colloquium on Group Theoretical Methods in Physics, Paris,
July 2002. To appear in the Proceedings of the Colloquium.}

\author{Edward Frenkel}

\address{Department of Mathematics, University of California,
  Berkeley, CA 94720, USA}


\maketitle



Representation theory of affine Kac-Moody algebras at the critical
level contains many intricate structures, in particular, the
hamiltonian structures of the KdV and modified KdV hierarchies and the
Miura transformation between them. In this talk I will describe these
structures and their deformations which will lead us to the deformed
Virasoro and $\W$--algebras and the integrable hierarchies associated
to them. I will also discuss briefly the relation of these matters to
the geometric Langlands correspondence.

It is a great honor for me to give this talk as the first recipient of
the Hermann Weyl Prize. Weyl was a pioneer of applications of symmetry
in quantum physics, a scientist who truly appreciated the beauty of
mathematics. He once said: {\em My work has always tried to unite the
true with the beautiful and when I had to choose one or the other, I
usually chose the beautiful.}

\section{Affine Kac-Moody algebras}

Let $\g$ be a simple Lie algebra over $\C$. Fix an invariant inner
product $\kappa$ on $\g$ and let $\ghat_\ka$ denote the
one-dimensional central extension of $\g \otimes \C((t))$ with the
commutation relations
\begin{equation}    \label{KM rel}
[A \otimes f(t),B \otimes g(t)] = [A,B] \otimes f(t) g(t) - (\kappa(A,B)
\on{Res} f dg) K,
\end{equation}
where $K$ is the central element. The Lie algebra $\ghat_\ka$ is the
{\em affine Kac-Moody algebra} associated to $\g$.

A representation of $\ghat_\ka$ on a complex vector space $V$ is
called {\em smooth} if for any vector $v \in V$ there exists $N \in
\Z_+$ such that $\g \otimes t^N \C[[t]] \cdot v = 0$. We also require
that $K$ acts on $V$ as the identity (since the space of invariant
inner products on $\g$ is one-dimensional, this is equivalent to a
more traditional approach whereby one fixes $\ka$ but allows $K$ to
act as the identity times a scalar, called the level).

Let $U_{\ka}(\ghat)$ be the quotient of the universal enveloping
algebra $U(\ghat_{\ka})$ of $\ghat_{\ka}$ by the ideal generated by
$(K-1)$. Define its completion $\wt{U}_{\ka}(\ghat)$ as follows:
$$
\wt{U}_{\ka}(\ghat) = \lim_{\longleftarrow} \;
U_{\ka}(\ghat)/U_{\ka}(\ghat) \cdot (\g \otimes t^N\C[[t]]).
$$
It is clear that $\wt{U}_{\ka}(\ghat)$ is a topological algebra which
acts on all smooth representations of $\ghat_{\ka}$, on which $K$ acts
as the identity. We shall recall the description of the center
$Z(\ghat)$ of $\wt{U}_{\ka_c}(\ghat)$ from \cite{FF:gd,F:new}.

First we need to introduce the notion of an oper.

\section{Opers}

We start with the definition of an $\sw_n$--oper. Let $X$ be a smooth
algebraic curve and $\Omega$ the line bundle of holomorphic
differentials on $X$. We will fix a square root $\Omega^{1/2}$ of
$\Omega$. An $\sw_n$--oper on $X$ is an $n$th order differential
operator acting from the holomorphic sections of $\Omega^{-(n-1)/2}$
to those of $\Omega^{(n+1)/2}$ whose principal symbol is equal to $1$
and subprincipal symbol is equal to $0$ (note that for these
conditions to be coordinate-independent, this operator must act from
$\Omega^{-(n-1)/2}$ to $\Omega^{(n+1)/2}$). Locally, we can choose a
coordinate $t$ and write this operator as
\begin{equation}    \label{sln-oper}
L = \partial_t^n+v_1(t) \partial_t^{n-2}+\ldots+v_{n-1}(t).
\end{equation}
It is not difficult to obtain the transformation formulas for the
coefficients $v_1(t),\ldots,$ $v_{n-1}(t)$ of this operator under
changes of coordinates.

For example, an $\sw_2$--oper is nothing but a Sturm-Liouville operator
of the form $\pa_t^2 + v(t)$ acting from $\Omega^{-1/2}$ to
$\Omega^{3/2}$. Under the change of coordinates $t = \varphi(s)$ we
have the following transformation formula
$$
v \mapsto \wt{v}, \qquad \wt{v}(s) = v(\varphi(s)) \left( \varphi'(s)
\right)^2 + \frac{1}{2} \{ \varphi,s \},
$$
where $$\{ \varphi,s \} = \frac{\varphi'''}{\varphi'} - \frac{3}{2}
\left( \frac{\varphi''}{\varphi'} \right)^2$$ is the Schwarzian
derivative. Operators of this form are known as {\em projective
  connections} (see, e.g., \cite{FB}, Sect. 7.2).

Drinfeld and V. Sokolov \cite{DS} have introduced an analogue of
operators \eqref{sln-oper} for a general simple Lie algebra
$\g$. Their idea was to replace the operator \eqref{sln-oper} by the
first order matrix differential operator
\begin{equation}    \label{sln-oper1}
\partial_t + \left( \begin{array}{ccccc}
0&v_1&v_2&\cdots&v_{n-1}\\
-1&0&0&\cdots&0\\
0&-1&0&\cdots&0\\
\vdots&\ddots&\ddots&\cdots&\vdots\\
0&0&\cdots&-1&0
\end{array}\right).
\end{equation}
Now consider the space of more general operators of the form
\begin{equation}    \label{sln-oper2}
\pa_t + \left( \begin{array}{ccccc}
*&*&*&\cdots&*\\
-1&*&*&\cdots&*\\
0&-1&*&\cdots&*\\
\vdots&\ddots&\ddots&\ddots&\vdots\\ 
0&0&\cdots&-1&*
\end{array} \right)
\end{equation}
The group of upper triangular matrices with $1$'s on the diagonal
acts on this space by gauge transformations
$$
\pa_t + A(t) \mapsto \pa_t + g A(t) g^{-1} - g^{-1} \pa_t g.
$$

It is not difficult to show that this action is free and each orbit
contains a unique operator of the form \eqref{sln-oper1}. Therefore,
locally, over a sufficiently small neighborhood $U$ of $X$, equipped
with a coordinate $t$, the space of $\sw_n$--opers on $U$ may be
identified with the quotient of the space of operators of the form
\eqref{sln-oper2} by the gauge action of the group of upper triangular
matrices.

Now let us generalize the latter definition to the case of an
arbitrary simple Lie algebra $\g$. Let $\g = \n_+ \oplus \h \oplus
\n_-$ be the Cartan decomposition of $\g$ and $e_i, h_i$ and $f_i,
i=1,\ldots,\ell$, be the Chevalley generators of $\n_+, \h$ and
$\n_-$, respectively. Then the analogue of the space of operators of
the form \eqref{sln-oper2} is the space of operators
\begin{equation}    \label{g-oper}
\pa_t - \sum_{i=1}^\ell f_i + {\mb v}(t), \qquad {\mb v}(t) \in \bb_+,
\end{equation}
where $\bb_+ = \h \oplus \n_+$. This space is preserved by the action
of the group of $N_+$--valued gauge transformations, where $N_+$ is
the Lie group corresponding to $\n_+$. Following \cite{DS}, we define
the space of $\g$--opers over a sufficiently small neighborhood $U$ of
$X$, equipped with a coordinate $t$, as the quotient of the space of
all operators of the form \eqref{g-oper} by the $N_+$--valued gauge
transformations. One shows that these gauge transformations act
freely, and one can find canonical representatives of each orbit
labeled by $\ell$--tuples of functions $v_1(t),\ldots,v_\ell(t)$. The
first of them, $v_1(t)$, transforms as a projective connection, and
$v_i(t), i>1$, transforms as a $k_i$--differential (i.e., a section of
$\Omega^{k_i}$), where $\{ k_i \}_{i=1,\ldots,\ell}$ are the orders of
the Casimir elements of $\g$, i.e., the generators of the center of
$U(\g)$. Using these transformations, one glues together the spaces of
$\g$--opers on different open subsets of $X$ and thus obtains the
notion of a $\g$--oper on $X$. A. Beilinson and V. Drinfeld have given
a more conceptual definition of opers on $X$ as $G$--bundles with a
connection and a reduction to the Borel subgroup $B_+$ (the Lie group
of $\bb_+$) satisfying a certain transversality condition, see
\cite{BD}, Sect. 3.

Denote by $\on{Op}_\g(D^\times)$ the space of $\g$--opers on the
(formal) punctured disc. This is the quotient of the space of
operators of the form \eqref{g-oper}, where ${\mb v}(t) \in
\bb_+((t))$, by the gauge action of $N_+((t))$.

Drinfeld and Sokolov have obtained $\on{Op}_\g(D^\times)$ as the
result of the hamiltonian reduction of the space of all operators of
the form $\pa_t + A(t), A(t) \in \g((t))$. The latter space may be
identified with a hyperplane in the dual space to the affine Lie
algebra $\ghat_{\nu_0}, \nu_0 \neq 0$, which consists of all linear
functionals taking value $1$ on $K$. It carries the Kirillov-Kostant
Poisson structure; the non-zero invariant inner product $\nu_0$ on
$\g$ appears as a parameter of this Poisson structure.

Applying the Drinfeld-Sokolov reduction, we obtain a Poisson structure
on the algebra $\on{Fun}(\on{Op}_\g(D^\times))$ of functions on
$\on{Op}_\g(D^\times)$. This Poisson algebra is called the {\em
classical $\W$--algebra} associated to $\g$. In the case when
$\g=\sw_n$, this Poisson structure is the (second) {\em
Adler--Gelfand--Dickey} Poisson structure. Actually, it is a member of
a two-dimensional family of Poisson structures on
$\on{Op}_{\sw_n}(D^\times)$ with respect to which the flows of the
$n$th KdV hierarchy are hamiltonian. We recall that in terms of the
first description of $\sw_n$--opers, i.e., as operators $L$ of the
form \eqref{sln-oper}, the KdV equations may be written in the Lax
form
\begin{equation}    \label{n-kdv}
\pa_{t_m} L = [(L^{m/n})_+,L], \qquad m>0, m \not{|} \; n,
\end{equation}
where $L^{1/n} = \pa_t + \ldots $ is the pseudodifferential operator
obtained by extracting the $n$th root of $L$ and $+$ indicates the
differential part of a pseudodifferential operator.

Drinfeld and Sokolov have defined an analogue of the KdV hierarchy on
the space of $\g$--opers for an arbitrary $\g$. The equations of this
hierarchy are hamiltonian with respect to the above Poisson structure
(in fact, they are bihamiltonian, but we will not discuss here the
other hamiltonian structure).

\section{The center}

Let us go back to the completed universal enveloping algebra
$\wt{U}_{\ka}(\ghat)$ and its center $Z_{\ka}(\ghat)$. Note that
$Z_{\ka}(\ghat)$ is a Poisson algebra. Indeed, choosing a non-zero
invariant inner product $\ka_0$, we may write a one-parameter
deformation of $\ka$ as $\ka + \ep \ka_0$. Given two elements, $A, B
\in Z_{\ka}(\ghat)$, we consider their arbitrary $\ep$--deformations,
$A(\ep), B(\ep) \in \wt{U}_{\ka+\ep \ka_0}(\ghat)$. Then the
$\ep$--expansion of the commutator $[A(\ep),B(\ep)]$ will not have a
constant term, and its $\ep$--linear term, specialized at $\ep=0$,
will again be in $Z_{\ka}(\ghat)$ and will be independent of the
deformations of $A$ and $B$. Thus, we obtain a bilinear operation on
$Z_{\ka}(\ghat)$, and one checks that it satisfies all properties of a
Poisson bracket.

Now we can describe the Poisson algebra $Z_{\ka}(\ghat)$. For a simple
Lie algebra $\g$ we denote by $^L \g$ its Langlands dual Lie algebra,
whose Cartan matrix is the transpose of that of $\g$ (note that this
duality only affects the Lie algebras of series $B$ and $C$, which get
interchanged). We have a canonical identification $^L \h = \h^*$.

Let $\ka_c$ be the critical inner product on $\g$ defined by the
formula $$\ka_c(x,y) = -\frac{1}{2} \on{Tr}_\g \on{ad} x \on{ad} y.$$
In the standard normalization of \cite{Kac}, the modules over
$\wt{U}_{\ka_c}(\ghat)$ on which $K$ acts as the identity are the
$\ghat$--modules of {\em critical level} $-h^\vee$, where $h^\vee$ is
the dual Coxeter number.

\begin{thm}[\cite{FF:gd,F:new}]\hfill    \label{center}

\noindent (1) {\em If $\ka \neq \ka_c$, then $Z_{\ka}(\ghat) = \C$.}

\noindent (2) The center {\em $Z_{\ka_c}(\ghat)$ is isomorphic, as a
Poisson algebra, to the classical $\W$--algebra $\on{Fun}({\on
Op}_{^L\g}(D^\times))$.}
\end{thm}

Thus, we recover the (second) Poisson structure of the $^L \g$--KdV
hierarchy from the center of the completed universal enveloping
algebra $\wt{U}_{\ka_c}(\ghat)$. Note that the two Poisson structures
appearing in the theorem depend on parameters: the inner products
$\kappa_0$ on $\g$ and $\nu_0$ on $^L \g$. In the above isomorphism
they have to agree in the obvious sense, namely, that the restriction
of $\ka_0$ to $\h$ is dual to the restriction of $\nu_0$ to $^L \h =
\h^*$.

For example, $\on{Op}_{\sw_2}(D^\times) = \{ \pa_t^2 - v(t) \}$, where
$v(t) = \sum_{n \in \Z} v_n t^{-n-2}$ is a formal Laurent
series. Therefore $\on{Fun}(\on{Op}_{\sw_2}(D^\times))$ is a
completion of the polynomial algebra $\C[v_n]_{n \in \Z}$. The Poisson
structure is that of the classical Virasoro algebra; it is uniquely
determined by the Poisson brackets between the generators
$$
\{ v_n,v_m \} = (n-m) v_{n+m} - \frac{1}{2} (n^3-n) \delta_{n,-m}.
$$

Under the above isomorphism, the generators $v_n$ are mapped to the
{\em Segal-Sugawara operators} $S_n$. Those are defined (for an
  arbitrary $\g$) by the formula
$$
S(z) = \sum_{n \in \Z} S_n \; z^{-n-2} = \sum_a {\mathbf :} J^a(z)^2
{\mathbf :},
$$
where
$$
J^a(z) = \sum_{n \in \Z} J^a_n \; z^{-n-1}, \qquad J^a_n = J^a \otimes
t^n,
$$
and $\{ J^a \}$ is an orthonormal basis of $\g$ with respect to
$\ka_0$.

For $\g=\sw_2$, the center $Z_{\ka_c}(\su)$ is a completion of the
polynomial algebra generated by $S_n, n \in \Z$. For a general $\g$,
we also have $\ell-1$ ``higher'' Segal--Sugawara operators $S^{(i)}_n,
i=2,\ldots,\ell, n \in \Z$, of orders equal to the orders of the
Casimirs of $\g$, and the center $Z_{\ka_c}(\ghat)$ is a completion of
the algebra of polynomials in these operators. However, explicit
formulas for $S^{(i)}_n$ with $i>1$ are unknown in general.

\section{Miura transformation}

In addition to operators of the form \eqref{sln-oper1}, it is useful
to consider the operators
\begin{equation}    \label{sln-moper}
\partial_t + \left( \begin{array}{ccccc}
u_1&0&0&\cdots&0\\
-1&u_2&0&\cdots&0\\
0&-1&u_3&\cdots&0\\
\vdots&\ddots&\ddots&\cdots&\vdots\\
0&0&\cdots&-1&u_n
\end{array}\right), \qquad \sum_{i=1}^n u_i = 0.
\end{equation}
It is easy to see that the operator \eqref{sln-moper} defines the same
oper as the operator \eqref{sln-oper1} (i.e., that they are gauge
equivalent under the action of the group of upper triangular matrices)
if and only if we have the following identity:
\begin{equation}    \label{split n}
\partial_t^n+v_1(t) \partial_t^{n-2}+\ldots+v_{n-1}(t) = (\pa_t+u_1(t))
\ldots (\pa_t+u_n(t)),
\end{equation}
This equation expresses $v_1,\ldots,v_{n-1}$ as differential
polynomials in $u_1,\ldots,u_\ell$. For example, for $n=2$ we have
\begin{equation}    \label{miura}
\pa_t^2 - v = (\pa_t - u)(\pa_t + u), \qquad \on{i.e.}, \qquad v = u^2
- u'.
\end{equation}
The latter formula is known as the {\em Miura
transformation}. R. Miura had found that this formula relates
solutions of the KdV equation to solutions of another soliton
equation, called the modified KdV, or mKdV, equation. This suggested
that the KdV equation has something to do with the second order
operators $\pa_t^2 - v$, because this formula appears in the splitting
of this operator into two operators of order one. It was this
observation that has subsequently led to the discovery of the inverse
scattering method (in a subsequent work by Gardner, Green, Kruskal and
Miura).

The Miura transformation may be viewed as a map from the space of the
first order operators $\{ \pa_t + u(t) \}$ to the space
$\on{Op}_{\sw_2}(D^\times) = \{ \pa_t^2 - v(t) \}$ of $\sw_2$--opers
(or projective connections) on $D^\times$. In order to make it
coordinate-independent, we must view the operator $\pa_t + u(t)$ as
acting from $\Omega^{-1/2}$ to $\Omega^{1/2}$, i.e., consider it as a
connection on the line bundle $\Omega^{-1/2}$. Denote the space of
such connections on $D^\times$ by $\on{Conn}_{\sw_2}(D^\times)$. Then
this map is actually a Poisson map $\on{Conn}_{\sw_2}(D^\times) \to
\on{Op}_{\sw_2}(D^\times)$ if we introduce the Poisson structure
$\on{Conn}_{\sw_2}(D^\times)$ by the formula
\begin{equation}    \label{hp}
\{ v_n,v_m \} = \frac{1}{2} n \delta_{n,-m},
\end{equation}
where $u(t) = \sum_{n \in \Z} u_n t^{-n-1}$. Thus, the algebra of
  functions on $\on{Conn}_{\sw_2}(D^\times)$, which is a completion of
  the polynomial algebra $\C[u_n]_{n \in \Z}$, is a Heisenberg--Poisson
  algebra.

Drinfeld and Sokolov have defined an analogue of the Miura
trans\-formation for an arbitrary simple Lie algebra $\g$. The role of
the operator $\pa_t+u(t)$ is now played by the operator $\pa_t + {\mb
u}(t)$, where ${\mb u}(t)$ takes values in $\h((t))$, considered as a
connection on the $H$--bundle $\Omega^{\rho^\vee}$. In other words,
under the change of variables $t = \varphi(s)$ it transforms as
follows:
$$
u \mapsto \wt{u}, \qquad \wt{u}(s) = u(\varphi(s)) \varphi'(s) -
\rho^\vee \left( \frac{\varphi''(s)}{\varphi'(s)} \right).
$$
Denote the space of such operators on the punctured disc by
$\on{Conn}_{\g}(D^\times)$. Then we have a natural map
$$
\mu: \on{Conn}_{\g}(D^\times) \to \on{Op}_\g(D^\times)
$$
which sends $\pa_t + {\mb u}(t)$ to the oper which is the gauge class
of the operator $\pa_t - \sum_{i=1}^\ell f_i + {\mb u}(t)$. This
is the Miura transformation corresponding to the Lie algebra $\g$.

Define a Poisson structure on $\on{Conn}_{\g}(D^\times)$ as
follows. Write $u_i(t) = \langle \al_i,{\mb u}(t) \rangle$ and $u_i(t)
= \sum_{n \in \Z} u_{i,n} t^{-n-1}$. Then set
$$
\{ u_{i,n},u_{j,m} \} = n \nu^{-1}_0(\al_i,\al_j) \delta_{n,-m},
$$
where $\nu_0^{-1}$ is the inner product on $\h^*$ induced by
$\nu_0|_{\h}$. Then the Miura transformation $\on{Conn}_{\g}(D^\times)
\to \on{Op}_\g(D^\times)$ is a Poisson map (if we take the Poisson
structure on $\on{Op}_\g(D^\times)$ corresponding to $\nu_0$).

Drinfeld and Sokolov \cite{DS} have defined the modified KdV hierarchy
corresponding to $\g$ on $\on{Conn}_{\g}(D^\times)$. The equations of
this $\g$--mKdV hierarchy are hamiltonian with respect to the above
Poisson structure. The Miura transformation intertwines the 
$\g$--mKdV and $\g$--KdV hierarchies.

\section{Wakimoto modules}

In \thmref{center} we identified the Poisson structure of the
$^L\g$--KdV hierarchy with the Poisson structure on the center
$Z_{\ka_c}(\ghat)$. It is natural to ask whether one can interpret in
a similar way the Poisson structure of the $^L\g$--mKdV hierarchy and
the Miura transformation. This can indeed be done using the Wakimoto
modules of critical level constructed in
\cite{Wak,FF:usp,FF:si,F:new}.

Let us briefly explain the idea of the construction of the Wakimoto
modules (see \cite{FB}, Ch. 10-11, and \cite{F:new} for more details).
Set $\bb_- = \h \oplus \n_-$. Given a linear functional $\chi: \h((t))
\to \C$, we extend it trivially to $\n_-((t))$ and obtain a linear
functional on $\bb_-((t))$, also denoted by $\chi$. Let $\C_\chi$ be
the corresponding one-dimensional representation of $\bb_-((t))$. We
would like to associate to it a smooth representation of $\ghat$. It
is clear that the induced module $\on{Ind}_{\bb_-((t))}^{\g((t))}
\C_\chi$ is not smooth. Therefore we need to modify the construction
of induction corresponding to a different choice of vacuum. In the
induced module the vacuum is annihilated by the Lie subalgebra
$\n_-((t))$, while in the Wakimoto module obtained by the
``semi-infinite'' induction the vacuum is annihilated by $t\g[t]
\oplus \n_+$.

However, when one applies the ``semi-infinite'' induction procedure
one has to deal with certain ``quantum corrections''. The effect of
these corrections is two-fold: first of all, the resulting module is a
module over the central extension of $\g((t))$, i.e., the affine
algebra $\ghat$, of critical level. Second, the parameters of the
module no longer behave as linear functionals on $\h((t))$, or,
equivalently, as elements of the space $\h^*((t)) dt = {}^L\h((t))dt$
of $\h^*$--valued one-forms on $D^\times$, but as connections on the
$^L H$--bundle $\Omega^\rho$. They are precisely the elements of the
space $\on{Conn}_{^L \g}(D^\times)$ which is a principal homogeneous
space over ${}^L\h((t))dt$.

Thus we obtain a family of smooth representations of
$\wt{U}_{\ka_c}(\ghat)$ (on which the central element $K$ acts as the
identity) parameterized by points of $\on{Conn}_{^L \g}(D^\times)$.
These are the {\em Wakimoto modules} of critical level. For $\chi \in
\on{Conn}_{^L \g}(D^\times)$ we denote the corresponding module by
$W_\chi$.

\smallskip

\noindent{\bf Example.} Let $\g=\sw_2$ with the standard basis $\{
e,h,f\}$. Consider the Weyl algebra with generators $a_n, a^*_n, n \in
\Z$, and relations $[a_n,a^*_m] = \delta_{n,-m}$. Let $M$ be the Fock
representation generated by a vector $|0\rangle$ such that $a_n
|0\rangle = 0, n\geq 0$ and $a^*_n |0\rangle = 0, n>0$. Set $a(z) =
\sum_{n \in \Z} a_n z^{-n-1}$, etc. Then for any Laurent series $u(t)$
the formulas
\begin{align*}
e(z) &= a(z), \\
h(z) &= - 2 {\mb :} a(z) a^*(z) {\mb :} + u(z),\\
f(z) &= - {\mb :} a(z) a^*(z)^2 {\mb :} + u(z) a^*(z) - 2 \pa_z a^*(z)
\end{align*}
define an $\su$--module structure on $M$. This is the Wakimoto module
attached to $\pa_t - u(t)$. One checks easily that in order for the
$h_n$'s to transform as the functions $t^n$ (or, equivalently, for
$h(z)dz$ to transform as a one-form), $\pa_t - u(t)$ needs to
transform as a connection on $\Omega^{-1/2}$.

One also checks that the Segal-Sugawara operator $S(z)$ acts on this
module as $u^2 - u'$, i.e., through the Miura transformation. This
statement has the following generalization for an arbitrary $\g$.

\begin{thm}[\cite{FF:gd,F:new}]
{\em The center $Z_{\ka_c}(\ghat)$ acts on $W_\chi, \chi \in
\on{Conn}_{^L \g}(D^\times)$, according to a character. The
corresponding point in $\on{Spec} Z_{\ka_c}(\ghat) = \on{Op}_{^L
\g}(D^\times)$ is $\mu(\chi)$, where $\mu: \on{Conn}_{^L \g}(D^\times)
\to \on{Op}_{^L \g}(D^\times)$ is the Miura transformation.}
\end{thm}

Thus, we obtain an interpretation of the Miura transformation as an
affine analogue of the Harish-Chandra homomorphism $Z(\g) \to
\C[\h^*]^W$. We remark that the map $\on{Conn}_{^L \g}(D^\times) \to
\on{Spec} Z_{\ka_c}(\ghat)$ that we obtain this way is manifestly
Poisson because the Wakimoto modules may be deformed away from the
critical level (this deformation gives rise to a Poisson structure on
$\on{Conn}_{^L \g}(D^\times)$ which coincides with the one introduced
above).

In summary, we have now described the phase spaces of both the
generalized KdV and mKdV hierarchies, their Poisson structures and a
map between them in terms of representations of affine Kac-Moody
algebras of critical level.

\section{Local Langlands correspondence for affine algebras}

The classical local Langlands correspondence aims to describe the
isomorphism classes of smooth representations of $G(F)$, where $G$ is a
reductive algebraic group and $F = {\mathbb Q}_p$ or ${\mathbb
F}_q((t))$, in terms of homomorphisms from the Galois group of $F$ to
the Langlands dual group $^L G$ (this is the group for which the sets
of characters and cocharacters of the maximal torus are those of $G$,
interchanged).

Let us replace ${\mathbb F}_q((t))$ by $\C((t))$ and $G$ by its Lie
algebra $\g$. Then we try to describe smooth representations of the
central extension of the loop algebra $\g((t))$ in terms of some
Galois data. But in the geometric context the Galois group should be
thought of as a sort of fundamental group. Hence we replace the notion
of a Galois representation by the notion of a $^L G$--local system, or
equivalently, a $^L G$--bundle with connection on $D^\times$.

The local Langlands correspondence in this context should be a
statement that to each $^L G$--bundle with connection on $D^\times$
corresponds a category of $\ghat$--modules. Here is an example of such
a statement in which the Wakimoto modules of critical level introduced
in the previous section play an important role (this is part of an
ongoing joint project with D. Gaitsgory).

First we introduce the notion of {\em nilpotent opers}. Those are
roughly those opers on $D^\times$ which have regular singularity at
the origin and unipotent monodromy around $0$. We denote the space of
nilpotent $\g$--opers by $\on{nOp}_\g$. For $\g=\sw_2$, its points are
the projective connections of the form $\pa_z^2 - v(z)$, where $v(t) =
\sum_{n\leq -1} v_n t^{-n-2}$. We have a residue map $\on{Res}:
\on{nOp}_\g \to \n$ which for $\g=\sw_2$ takes the form $\pa_t^2 -
v(t) \mapsto v_{-1}$.

Recall that for a nilpotent element $x \in {}^L \g$, the {\em Springer
fiber} of $x$ is the variety of all Borel subalgebras of $^L \g$
containing $x$. For example, the Springer fiber at $0$ is just the
flag variety of $^L \g$.

\begin{lem}
The set of points of the fiber $\mu^{-1}(\rho)$ of the Miura
transformation over a nilpotent $^L \g$--oper $\rho$ is in bijection
with the set of points of the Springer fiber of $\on{Res}(\rho) \in
{}^L\n$.
\end{lem}

For example, the set of points of the fiber $\mu^{-1}(\rho)$ of the
Miura transformation over a regular $^L \g$--oper is the set of points
of the flag variety of $^L \g$.

Now fix $\rho \in \on{nOp}_{^L \g}$ and consider the category ${\mc
C}_\rho$ of $\ghat$--modules of critical level on which
$Z_{\ka_c}(\ghat) \simeq \on{Fun}(\on{Op}_{^L\g}(D^\times))$ acts
through the central character $Z_{\ka_c}(\ghat) \to \C$ corresponding
to $\rho$ and such that the Lie subalgebra $(t\g[[t]] \oplus \n_+)
\subset \ghat$ acts locally nilpotently and the Lie subalgebra $\h$
acts with integral generalized eigenvalues. Then, according to a
conjecture of Gaitsgory and myself, the derived category of ${\mc
C}_\rho$ is equivalent to the derived category of quasicoherent
sheaves on the Springer fiber of $\on{Res}(\rho)$ (more precisely, the
corresponding DG-scheme).

In particular, under this equivalence the skyscraper sheaf at a point
of the Springer fiber of $\on{Res}(\rho)$, which is the same as a
point $\chi$ of $\on{Conn}_{^L\g}(D^\times)$ projecting onto $\rho$
under the Miura transformation, should correspond to the Wakimoto
module $W_\chi$. Thus, the above conjecture means that, loosely
speaking, any object of the category ${\mc C}_\rho$ may be
``decomposed'' into a ``direct integral'' of Wakimoto modules.

\section{A $q$--deformation}

Now we wish to define $q$--deformations of the structures described in
the previous sections. In particular, we wish to introduce
$q$--analogues of opers and connections (together with their Poisson
structures) and of the Miura transformation between them. We also wish
to define $q$--analogues of the KdV hierarchies. For that we replace
the universal enveloping algebra of the affine Lie algebra $\ghat$ by
the corresponding {\em quantized enveloping algebra} $U_q(\ghat)$.

Then the center of $U_q(\ghat)$ at the critical level (with its
Poisson structure defined in the same way as in the undeformed case)
should be viewed as $q$--analogue of the algebra of functions on opers
(i.e., the classical $\W$--algebra). One the other hand, parameters of
Wakimoto modules should be viewed as $q$--analogues of connections,
and the action of the center of $U_q(\ghat)$ on Wakimoto modules
should give us a $q$--analogue of the Miura transformation.

In \cite{FR1}, N. Reshetikhin and I have computed these structures in
the case when $\g=\sw_n$ (we used the Wakimoto modules over
$U_q(\wh{\sw}_n)$ constructed in \cite{AOS}). Let us describe the
results. The $q$--analogues of $\sw_n$--opers are $q$--difference
operators of the form
\begin{equation}    \label{qsln-oper}
L_q = D^n + t_1(z) D^{n-1} + \ldots + t_{n-1}(z) D + 1,
\end{equation}
where $(Df)(z) = f(zq^2)$. The $q$--analogues of connections are
operators $D + {\mb \Lambda}(z)$, where ${\mb \Lambda} =
(\Lambda_1,\ldots,\Lambda_n)$ and $\prod_{i=1}^n \Lambda_i(z) =
1$. The $q$--analogue of the Miura transformation is the formula
expressing the splitting of the operator \eqref{qsln-oper} into a
product of first order operators
$$
D^n + t_1(z) D^{n-1} + \ldots + t_{n-1}(z) D + 1 = (D+\Lambda_1(z))
\ldots (D+\Lambda_n(z)).
$$

For example, for $\g=\sw_2$ the $q$--Miura transformation is
\begin{equation}    \label{baxter}
t(z) = \Lambda(z) + \Lambda(zq^2)^{-1}.
\end{equation}
Note that in the limit $q = e^h, h \to 0$ we have $t(z) = 2 + 4 h^2
v(z) z^2 + \ldots$ and $\Lambda(z) = e^{2 h u(z) z}$ so that we obtain
the ordinary Miura transformation $v = u^2 - u'$.

The Poisson structures with respect to which this map is Poisson are
given by the formulas \cite{FR1}
$$
\{ \Lambda(z),\Lambda(w) \} = (q-q^{-1}) f(w/z) \; \Lambda(z)
\Lambda(w), \qquad f(x) = \sum_{n\in\Z} \frac{q^n - q^{-n}}{q^n +
  q^{-n}} \; x^n,
$$
$$
\{ t(z),t(w) \} = (q-q^{-1}) \left( f(w/z) \; t(z)t(w) + \delta\left(
\frac{w}{zq^2} \right) - \delta\left( \frac{wq^2}{z} \right) \right),
$$
where $\delta(x) = \sum_{n \in \Z} x^n$. Analogous formulas for $\sw_n,
n>2$, may be found in \cite{FR1}.

The equations of the $q$--analogue of the $n$th KdV hierarchy are
given by the formulas
$$
\pa_{t_m} L_q = \left[ (L_q^{m/n})_+,L_q \right], \quad \quad
m>0, m \not{|} \; n.
$$
These equations are hamiltonian with respect to the above Poisson
structure \cite{F:qkdv}.

For simple Lie algebras other than $\sw_n$, Wakimoto modules have not
yet been constructed. Nevertheless, in \cite{FR2} we have generalized
the above formulas to the case of an arbitrary $\g$.

Another approach is to define $q$--analogues of the classical
$\W$--algebras by means of a $q$--analogue of the Drinfeld-Sokolov
reduction. This has been done in \cite{FRS,SS}.

If we write $\Lambda(z) = Q(zq^{-2})/Q(z)$, then formula
\eqref{baxter} becomes
$$
t(z) = \frac{Q(zq^2)}{Q(z)} + \frac{Q(zq^{-2})}{Q(z)},
$$
which is the Baxter formula for the Bethe Ansatz eigenvalues of the
transfer-matrix in the XXZ spin chain model. This is not
coincidental. The transfer-matrices of spin models may be obtained
from central elements of $U_q(\ghat)$ of critical level and the Bethe
eigenvectors can be constructed using Wakimoto modules. Then the
formula for the eigenvalues becomes precisely the formula for the
$q$--Miura transformation (see \cite{FR2}). In the quasi-classical
limit we obtain that the formula for the Bethe Ansatz eigenvalues of
the Gaudin model (the quasi-classical limit of the XXZ model) may be
expressed via the (ordinary) Miura transformation formula. This was
explained (in the case of an arbitrary simple Lie algebra $\g$) in
\cite{FFR,F:icmp}.

\section{The $q$--characters}

Let $\on{Rep} \, U_q(\ghat)$ be the Grothendieck ring of
finite-dimensional representations of $U_q(\ghat)$.  N.Reshetikhin and
M.Semenov-Tian-Shansky \cite{RS} have given an explicit construction
of central elements of $U_q(\ghat)$ of critical level. It amounts to a
homomorphism from $\on{Rep} \, U_q(\ghat)$ to $Z_q((z))$, where $Z_q$
is the center of $U_q(\ghat)$ at the critical level. Combining this
homomorphism with the $q$--Miura transformation defined in \cite{FR2},
we obtain an injective homomorphism
$$
\chi_q: \on{Rep} \, U_q(\ghat) \to \Z[Y_{i,a}^{\pm 1}]_{i=1,\ldots,\ell;a
  \in \C^\times}
$$
which we call the $q$--character homomorphism (see \cite{FR3}). It
should be viewed as a $q$--analogue of the ordinary character
homomorphism
$$
\chi: \on{Rep} \, U(\g) \to \Z[y_i^{\pm 1}]_{i=1,\ldots,\ell},
$$ where the $y_i$'s are the fundamental coweights of $\g$. Under the
forgetful homomorphism $Y_{i,a} \mapsto y_i$, the $q$--character of a
$U_q(\ghat)$--module $V$ becomes the ordinary character of the
restriction of $V$ to $U_q(\g)$ (specialized at $q=1$). For instance,
if $V(a)$ is the two-dimensional representation of $U_q(\su)$ with the
evaluation parameter $a$, then $\chi_q(V(a)) = Y_a + Y_{aq^2}^{-1}$,
which under the forgetful homomorphism become the character $y+y^{-1}$
of the two-dimensional representation of $U_q(\sw_2)$ (or
$U(\sw_2)$). In \cite{FM2}, the notion of $q$--characters was extended
to the case when $q$ is a root of unity.

When $\g=\sw_n$, the variables $Y_{i,a}$ correspond to the series
$\Lambda_i(z)$ introduced above by the following rule: $\Lambda_i(za)
\mapsto Y_{i,aq^{-i+1}} Y_{i-1,aq^{-i+2}}^{-1}$, where $Y_0=Y_n=1$.

One gains a lot of insight into the structure of finite-dimensional
representations of quantum affine algebras by analyzing the
$q$--characters. For instance, it was conjectured in \cite{FR3} and
proved in \cite{FM1} that the image of the $q$--character homomorphism
is equal to the intersection of the kernels of certain {\em screening
operators}, which come from the hamiltonian interpretation of the
$q$--character homomorphism as the $q$--Miura transformation. This
enabled us to give an algorithm for the computation of the
$q$--characters of the fundamental representations of $U_q(\ghat)$
\cite{FM1}.

H. Nakajima \cite{N} has interpreted the $q$--characters in terms of
the cohomologies of certain quiver varieties. Using this
interpretation, he was able to describe the multiplicities of
irreducible representations of $U_q(\ghat)$ inside tensor products of
the fundamental representations.

\section{Deformed $\W$--algebras}

As we mentioned above, the $q$--character homomorphism expresses the
eigenvalues of the transfer-matrices in spin models obtained via the
Bethe Ansatz (see \cite{FR2}). The transfer-matrices form a
commutative algebra in a quantum object, the quantized enveloping
algebra $U_q(\ghat)$. But now we know that the algebra of
transfer-matrices carries a Poisson structure and that the
$q$--character homomorphism defines a Poisson map, i.e., the
$q$--Miura transformation. This immediately raises the question as to
whether the algebra of transfer-matrices (already a quantum, albeit
commutative, algebra) may be further quantized. This ``second
quantization'' was defined in \cite{SKAO,FF:w,AKOS,FR2}, and it leads
us to deformations of the $\W$--algebras.

The deformed $\W$--algebra $\W_{q,t}(\g)$ is a two-parameter family of
associative algebras. It becomes commutative in the limit $t \to 1$,
where it coincides with the $q$--deformed classical $\W$--algebra
discussed above. In another limit, when $t=q^\beta$ and $q \to 1$, one
obtains the conformal $\W$--algebras which first appeared in conformal
field theory (see \cite{FF:laws}).

For $\g=\sw_2$ we have the deformed Virasoro algebra
$\W_{q,t}(\sw_2)$. It has generators $T_n, n \in \Z$, satisfying the
relations
$$
f\left( \frac{w}{z} \right)  T(z) T(w) - f \left( \frac{z}{w} \right)
T(w) T(z) \\ = (q-q^{-1})(t-t^{-1}) \left( \delta \left(
\frac{w}{zq^2t^2} \right) - \delta \left( \frac{wq^2t^2}{z} \right)
\right)
$$
where $T(z) = \sum_{n \in \Z} T_n z^{-n}$ and
$$
f(z) = \frac{1}{1-z}
\frac{(zq^2;q^4t^4)_\infty (zt^2;q^4t^4)_\infty}{(zq^4t^2;q^4t^4)_\infty
  (zq^2 t^4;q^4t^4)_\infty}, \qquad (a;b)_\infty = \prod_{n=0}^\infty
(1-ab^n).
$$
There is also an analogue of the Miura transformation (free field
realization) given by the formula
$$
T(z) = {\mathbf :} \Lambda(z) {\mathbf :} + {\mathbf :} \Lambda(zq^2
t^2)^{-1} {\mathbf :} \; ,
$$
where $\Lambda(z)$ is the exponential of a generating function of
generators of a Heisenberg algebra.

Just like the Virasoro and other conformal $\W$--algebras, which are
symmetries of CFT, the deformed $\W$--algebras appear as dynamical
symmetry algebras of various models of statistical mechanics (see
\cite{LP,Miwa}).

\end{document}